\documentclass[12pt]{amsart}
\usepackage{color}
\usepackage[all]{xy}
\usepackage{amscd,amsmath,latexsym,amssymb}
\usepackage[mathcal]{euscript}
\usepackage{mathrsfs}

\date{August 15, 2018}

\calclayout
\newtheorem{dummy}{anything}[section]

\newtheorem*{thma}{Theorem A}

\theoremstyle{definition}%%Change Theoremstyle
\newtheorem{definition}[dummy]{Definition}

  \newtheorem{remark}[dummy]{Remark}

  \newtheorem*{acknowledgement}{Acknowledgement}

\newcommand{\sN}{\mathscr N}

\newcommand{\bZ}{\mathbb Z}

\newcommand{\bbR}{\mathbb R}

\newcommand{\cy}[1]{\bZ/{#1}}

\newcommand{\bd}{\partial}

\newcommand{\La}{\Lambda}
\newcommand{\mmatrix}[4]{\bigg (\hskip-4pt\vcenter
{\xymatrix@C-2pc@R-2pc{#1&#2\\#3&#4} }\hskip-2pt
\bigg )}

\begin{document}

\title[Orientable $4$-dimensional Poincar\'e Complexes are Reducible]
{Orientable $4$-dimensional Poincar\'e  Complexes have Reducible Spivak Fibrations}

\author{Ian Hambleton}
\address{\vbox{\hbox{Department of Mathematics \& Statistics, McMaster University}\hbox{Hamilton, Ontario L8S 4K1, Canada}}}

\email{hambleton@mcmaster.ca }

\thanks{Research partially supported by NSERC Discovery Grant A4000}

\begin{abstract}
We show that  the Spivak normal fibration of an orientable $4$-dimensional  Poincar\'e complex has a vector bundle reduction. \end{abstract}

\maketitle

\section{Introduction}
A Poincar\'e complex (\emph{PD-complex}), as introduced by Wall \cite[p.~214]{wall-pc1}, is a (connected) finitely dominated CW complex $X$ equipped with:

\begin{enumerate}
\item a homomorphism $w \colon \pi_1(X) \to \{\pm 1\}$ defining a twisted $\La := \bZ\pi_1(X)$  module structure $\bZ^t$ on $\bZ$.
\item an integer $n$ and a class $[X] \in H_n(X;\bZ^t)$ such that
\item for all integers $r\geq0$, cap product with $[X]$ induces an isomorphism
$$[X] \frown \colon H^r(X;\La) \to H_{n-r}(X; \La\otimes \bZ^t)\ .$$
\end{enumerate} 
The integer $n=\dim X$ is called the \emph{dimension} of $X$.
It follows from the foundational results of Kirby and Siebenmann \cite[Annex 3]{kirby-siebenmann1} that every closed topological $n$-manifold has the homotopy type of a Poincar\'e complex of dimension $n$ (see the discussion in Wall \cite[\S 17B]{wall-book}). In the manifold case, the homomorphism $w \colon \pi_1(X) \to \{\pm 1\}$ is given by the first Stiefel-Whitney class. Accordingly, a PD-complex $X$ is called \emph{orientable} if its homomorphism $w$ is trivial.

Spivak \cite{Spivak:1967} discovered that every simply-connected PD-complex $X$ with $\dim X = n$ has an associated spherical fibration, denoted $\nu_X$, which is unique up to stable fibre homotopy equivalence. It is constructed by embedding  $X$ in a high-dimensional Euclidean space $\bbR^{n+k}$ ($k \gg n$), and considering the fibration homotopic to the projection map $p\colon \bd N \to X$ from the boundary of a regular neighbourhood $N\subset \bbR^{n+k}$. The duality properties of $X$ imply that the fibres of $p$ are homotopy equivalent to $S^{k-1}$. The definition and the uniqueness statement were generalized by Wall \cite[\S 3]{wall-pc1} to all PD-complexes, and $\nu_X$ is now called the \emph{Spivak normal fibration} of  $X$.

 In the smooth manifold case, $\nu_X$ is the spherical fibration associated to the sphere bundle of  the (stable) normal $k$-vector bundle of $X$. For topological manifolds, the corresponding notion is the (stable) normal $\bbR^k$-bundle ($k \gg n$), and its sub-bundle with fibres $\bbR^k - \{0\} \simeq S^{k-1}$.

After the further development of geometric surgery theory, due to Browder, Milnor, Novikov, Sullivan and Wall, the normal structures on PD-spaces and manifolds were re-expressed via classifying spaces (see \cite[\S 10 and \S 17B]{wall-book}, \cite{kirby-siebenmann1}, \cite{madsen-milgram1}, \cite{kirby-taylor1}). One outcome was the construction of a sequence of classifying spaces
$$BO \to BPL \to BTOP \to BG$$
relating smooth, PL, and topoogical bundles to spherical fibrations. In particular, the (stable) Spivak normal fibre space $\nu_X$ is classified by a map $\nu_X \colon X \to BG$. 

\begin{definition} We say that  PD-complex $X$ has a \emph{reducible Spivak normal fibration} if the classifying map $\nu_X \colon X \to BG$ lifts to a map $\tilde\nu_X \colon X \to BTOP$. 
\end{definition}
Similarly, we say that the Spivak normal fibre space is reducible to a vector bundle if $\nu_X$ lifts to a map $\tilde\nu_X\colon X \to BO$. The lifting obstruction is given by the image of $\nu_X$ in 
$[X, B(G/TOP)]$ or $[X, B(G/O)]$, respectively. In dimensions $\geq 5$, these are different problems, but if $\dim X \leq 4$ these two obstruction groups are the same because
$$[X, B(G/O)] = [X, B(G/PL)] = [X, B(G/O)] \cong H^3(X;\cy 2), \quad \text{if\ } \dim X \leq 4.$$
This is explained in Kirby-Taylor \cite[\S 2]{kirby-taylor1}. In other words, the obstruction to reducibility for the Spivak normal fibration of a PD-complex $X$ in dimensions $\leq 4$ is a single characteristic class $k_3(X) \in H^3(X;\cy 2)$.

\begin{thma} Let $X$ be an Poincar\'e complex. If $\dim X \leq 3$, or $\dim X = 4$ and $X$ is orientable, then the Spivak normal fibration of $X$ is reducible to a vector bundle.\end{thma}

\begin{remark} The dimension $4$ case was known to the experts (see the statements in Spivak \cite[p.~95]{Spivak:1967} and Kirby-Taylor \cite[p.~10]{kirby-taylor1}), but  Land \cite{Land:2017} pointed out the  lack of a proof in the literature, and provided his own argument. 
For dimensions $\leq 2$ the result is immediate, and the dimension $3$ cases follow easily from the dimension $4$ statement. In general, non-oriented PD-complexes in dimensions $\geq 4$ do not have reducible Spivak normal fibrations (see Hambleton and Milgram \cite{hmilgram2} for explicit examples in every even dimension $\geq 4$). The first non-reducible \emph{orientable} example occurs in dimension 5 (see Gitler and Stasheff \cite{Gitler:1965}).
\end{remark}

\begin{acknowledgement} I would like to thank Wolfgang L\"uck for asking me about this result at a conference in Regensburg (September, 2017).  Andrew Ranicki and Larry Taylor later outlined alternate arguments, both different from the approach used by Markus Land, and  different from the proof provided in this note.
\end{acknowledgement}

\section{The proof of Theorem A}

Here is a short argument to show that an orientable 4-dimensional  Poincar\'e complex  has a reducible Spivak normal fibration. The proof is essentially contained in \cite{hmilgram2}.

\medskip
\noindent
1. Suppose that $X$ is an orientable 4-dimensional  PD-complex such that $\nu_X$ is not reducible. Then by Poincar\'e duality there is a class $e \in H^1(X.\cy 2)$ such that 
$$\langle k_3(X) \cup e , [X]\rangle  \neq 0,$$
where $k_3(X)$ denotes the pullback to $X$ of the first exotic characteristic class.
% (see \cite[Section 5]{hmilgram2}).

\medskip
\noindent
2. Let $f\colon X \to RP^{\infty}$ represent the cohomology class $e \in H^1(X;\cy 2)$. Then the element $ 0 \neq (X,f) \in \sN^{PD}_4(RP^\infty)$ has Arf invariant $A(X,f) = 1$ (see \cite{hmilgram2}, Corollary 4.2, Corollary 5.3, and Theorem 5.6).

\medskip
\noindent
3. By low-dimensional surgery, we may assume that $\pi_1(X) = \cy 2$ and that $f\colon X \to RP^{\infty}$ classifies its universal covering $\widetilde X \to X$ (see Wall \cite[Corollary 2.3.2]{wall-pc1}   to justify this much Poincar\'e surgery).

\medskip
\noindent
4. The form $B(a,b) = \langle a\,\cup\,T^*b, [X]\rangle $  is a symmetric unimodular bilinear form on $H^2(\widetilde X, \bZ)$, where $T$ denotes the non-trivial covering involution. The form $B$ is even (see Bredon \cite[ Chap VII, Theorem 7.4]{bredon1}).

\medskip
\noindent
5. The invariant $A(X,f)$ is the Arf invariant associated to the Browder-Livesay quadratic map $q$
 (see \cite[\S 4]{browder-livesay1}, and \cite[Theorem 1.4]{hmilgram2}), which refines the mod 2 reductions of $B$. By \cite[Lemma 4.6]{browder-livesay1},  we have
 $$q(a) \equiv \frac{B(a,a)}{2} \ {\rm (mod \ 2)}$$ 
 since $T\colon \widetilde X \to \widetilde X$ is orientation preserving.
 But $B$ is an even unimodular symmetic bilinear form, so  the Arf invariant obtained in this way is zero, and we have a contradiction. \qed

\begin{remark}
 To obtain the reducibility results for $3$-dimensional PD-complexes, one can make an appropriate circle bundle construction (which does not affect reducibility) resulting in orientable $4$-dimensional PD-complexes, and then apply Theorem A.
\end{remark}

%\bibliographystyle{ih}
%\bibliography{ihmain,4D}
%\end{document}
\providecommand{\bysame}{\leavevmode\hbox to3em{\hrulefill}\thinspace}
\providecommand{\MR}{\relax\ifhmode\unskip\space\fi MR }
% \MRhref is called by the amsart/book/proc definition of \MR.
\providecommand{\MRhref}[2]{%
  \href{http://www.ams.org/mathscinet-getitem?mr=#1}{#2}
}
\providecommand{\href}[2]{#2}

\end{document}